\documentstyle{article}
\begin{document}
\baselineskip+8pt
 \small \begin{center}{\textit{ In the name of
 Allah, the Beneficent, the Merciful}}\end{center}

\begin{center}{\bf  A radius of absolute convergence for power series in many variables }
\end{center}
\begin{center}{\bf  Ural Bekbaev  \\
 Turin Polytechnical University in Tashkent,\\ INSPEM, Universiti Putra Malaysia.\\
e-mail: bekbaev@science.upm.edu.my}
 \end{center}

 \begin{abstract}{In this paper for power series in many (real or complex)variables  a radius of (absolute) convergence
 is offered. This radius can be evaluated by a formula similar to Cauchy-Hadamard formula and in one
 variable case they are same.
 }\end{abstract}


  \vspace{0.5cm}

 In [1] a matrix representation of polynomial maps was offered and by the use of a new
 product of matrices the matrix representation of composition  of polynomial maps was given. In this paper an application of this product
 to power series in many (real or complex)variables is presented. Namely, by the use of this
 product power series in many variables are presented
 in the form of power series in one variable. Then by the use of $\rho$-norm (defined in [1]) and the
 Cauchy-Hadamard
 type formula a radius of
 (absolute) convergence for such series is introduced and investigated. In one variable case it is the same Cauchy-Hadamard
 formula for radius of
 convergence of power series in one variable.

 Here are some definitions and results related to the new product introduced in [1].

 For a positive integer $n$ let $I_n$ stand for all row $n$-tuples with nonnegative integer entries with the
 following linear order: $\beta=
 (\beta_1,\beta_2,...,\beta_n)<\alpha=(\alpha_1,\alpha_2,...,\alpha_n)$ if and only if $\vert \beta\vert < \vert \alpha\vert$ or
$\vert \beta\vert = \vert \alpha\vert$ and $\beta_1> \alpha_1$ or
$\vert \beta\vert = \vert \alpha\vert$, $\beta_1= \alpha_1$ and
$\beta_2> \alpha_2$ etcetera, where $\vert \alpha\vert$ stands for
$\alpha_1+\alpha_2+...+\alpha_n $.

It is clear that for $\alpha, \beta, \gamma \in I_n$ one has
$\alpha < \beta $ if and only if $\alpha + \gamma < \beta
+\gamma$. We write $\beta \ll \alpha$ if $\beta_i \leq \alpha_i$
for all $i=1,2,...,n$, $\left(\begin{array}{c}
  \alpha \\
  \beta \\
\end{array}\right)$ stands for $\frac{\alpha!}{\beta!
(\alpha -\beta)!}$, $\alpha !=\alpha_1!\alpha_2!...\alpha_n!$.

 In future $n$ and $n'$ are assumed to be any fixed positive
 integers. Let $F$ stand for the field of real or complex numbers.

 For any nonnegative
integer numbers $p,p'$ let  $M_{n,n'}(p,p';F)=M(p,p';F)$ stand for
all $"p\times p'"$ size matrices $A=(A_{\alpha,\alpha'})_{\vert
\alpha \vert=p, \vert \alpha'\vert=p'}$ ($\alpha$ presents row,
$\alpha'$ presents column and $\alpha\in I_n,\alpha'\in I_{n'}$)
with entries from $F$. Over such kind matrices in addition to the
ordinary sum and product of matrices we consider the following
"product" $\bigodot $ as well:

{\bf Definition 1.} If $A\in M(p,p';F)$ and $B\in M(q,q';F)$ then
$A\bigodot B=C\in M(p+q,p'+q';F)$ that for any
$\vert\alpha\vert=p+q$, $\vert\alpha'\vert=p'+q'$, where
$\alpha\in I_n,\alpha'\in I_{n'}$,
$$C_{\alpha,\alpha'}=\sum_{\beta,\beta'}\left(\begin{array}{c}
  \alpha \\
  \beta \\
\end{array}\right)
    A_{\beta,\beta'}B_{\alpha-\beta,\alpha'-\beta'}
$$, where the sum is taken over all $\beta\in I_n,\beta'\in I_{n'}$, for which $\vert
\beta\vert=p$, $\vert \beta'\vert=p'$, $\beta\ll \alpha$ and
$\beta'\ll \alpha'$.

Let us agree that $h$ ($H$, $v$, $V$) stands for any element of
$M(0,1;F)$ (respect. $M(0,p;F)$, $M(1,0;F)$, $M(p,0;F)$ , where
$p$ may be any nonnegative integer). We use $E_k$ for $"k\times
k"$ size ordinary unit matrix from $M_{n,n}(k,k;R)$. For the sake
of convenience it will be assumed that $A_{\alpha,\alpha'}=0$
($\alpha !=\infty$) whenever $\alpha \notin I_n$ or $\alpha'
\notin I_{n'}$ (respect. $\alpha \notin I_n$).

{\bf Proposition 1.} For the above defined product the following
are true.

1. $A\bigodot B=B\bigodot A$.

2. $(A+B)\bigodot C=A\bigodot C+ B\bigodot C$.

3. $(A\bigodot B)\bigodot C=A\bigodot (B\bigodot C)$

4. $ (\lambda A)\bigodot B=\lambda (A\bigodot B)$ for any
$\lambda\in F$

5. $A\bigodot B=0$ if and only
    if $A=0$ or $B=0$.

6. $A(B\bigodot H)=(AB)\bigodot H$

7. $(E_k\bigodot V)A=A\bigodot V$

In future $A^{(m)}$ means the $m$th power of the matrix $A$ with
respect to the new product.

{\bf Proposition 2.} If $h=(h_1,h_2,...,h_{n})\in M(0,1;F)$,
$v=(v_1,v_2,...,v_n)\in M(1,0;F)$, then
$$(h^{(m)})_{0,\alpha'}=\left(\begin{array}{c}
  m \\
  \alpha'\\
\end{array}\right)h^{\alpha'}, \hspace{1cm} (v^{(m)})_{\alpha,0}=m!v^{\alpha}$$, where $h^{\alpha}$ stands for
$h_1^{\alpha_1}h_2^{\alpha_2}...h_n^{\alpha_n}$

 In future let $\rho \geq 1$
 be any fixed real number and $\varrho$ stand for the real number for which $\frac{1}{\rho}+\frac{1}{\varrho}=1$.
We consider the following $\rho$-norm of elements $A\in
M(p,p';F)$:

{\bf Definition 2.}
$$\|A\|=\|A\|_{\rho}=(\sum_{\alpha,\alpha'}\frac{\vert A_{\alpha,\alpha'}\vert ^{\rho}}{\alpha!(p!p'!)^{\rho
-1}})^{1/\rho}$$

In the case of $\rho =\infty$ the $\rho$-norm is defined by
$$\|A\|=\|A\|_{\infty}=\frac{\sup_{\alpha, \alpha'}\vert
A_{\alpha,\alpha'}\vert }{p!p'!}$$

{\bf Theorem 1.} 1. If $A,B\in Mat(p,p';F)$ and $\lambda \in F$
then

a)$\|A\|=0$ if and only if $A=0$,

b)$\|\lambda A\| =|\lambda |\|A\|$,

c) $\|A+B\|\leq \|A\|+ \|B\|$.

2. For any nonnegative integer numbers $p$, $p'$, $q$ and $q'$
there is such a positive number $\lambda(p, p', q, q')$ that for
any $A\in Mat(p,p';F)$, $B\in Mat(q,q';F)$ the following
inequality is valid:
 $$\lambda(p, p', q, q')\|A\|\|B\| \leq \|A\bigodot B\|\leq \|A\|\|B\|$$

{\bf Proof.}  Here is a proof of part 2. First let us show the
inequality $\|A\bigodot B\|\leq \|A\|\|B\|$.

     Due to the H\"{o}lder inequality  for $A\bigodot B=C$ one has
$$\vert C_{\alpha,\alpha'}\vert =\vert\sum_{\beta \ll\alpha,\beta'\ll\alpha'}\left(\begin{array}{c}
  \alpha \\
  \beta \\
\end{array}\right)A_{\beta,\beta'}B_{\alpha-\beta,\alpha'-\beta'}\vert \leq
\alpha !\sum_{\beta \ll\alpha,\beta'\ll\alpha'}\frac{\vert
A_{\beta,\beta'}B_{\alpha-\beta,\alpha'-\beta'}\vert}{(\beta
!(\alpha-\beta)!)^{1/\rho}}\frac{1}{(\beta
!(\alpha-\beta)!)^{1/\varrho}} \leq $$ $$\alpha !(\sum_{\beta
\ll\alpha,\beta'\ll\alpha'}\frac{\vert
A_{\beta,\beta'}B_{\alpha-\beta,\alpha'-\beta'}\vert^{\rho}}{\beta
!(\alpha-\beta)!})^{1/\rho}(\sum_{\beta
\ll\alpha,\beta'\ll\alpha'}(\frac{1}{(\beta
!(\alpha-\beta)!)^{1/\varrho}})^{\varrho})^{1/\varrho}=$$
$$\alpha !(\sum_{\beta
\ll\alpha,\beta'\ll\alpha'}\frac{\vert
A_{\beta,\beta'}\vert^{\rho}}{\beta !}\frac{\vert
B_{\alpha-\beta,\alpha'-\beta'}\vert^{\rho}}{(\alpha-\beta)!})^{1/\rho}(\sum_{\beta
\ll\alpha}\frac{1}{\beta
!(\alpha-\beta)!}\sum_{\beta'\ll\alpha'}1)^{1/\varrho} \leq$$
$$\alpha !(\sum_{\beta
\ll\alpha,\beta'\ll\alpha'}\frac{\vert
A_{\beta,\beta'}\vert^{\rho}}{\beta !}\frac{\vert
B_{\alpha-\beta,\alpha'-\beta'}\vert^{\rho}}{(\alpha-\beta)!})^{1/\rho}(\left(
\begin{array}{c}
  p+q \\
  p \\
\end{array}\right)\frac{1}{\alpha !}\left(
\begin{array}{c}
  p'+q' \\
  p' \\
\end{array}\right))^{1/\varrho}$$ as far as according to
Proposition 1 one has $\sum_{\beta \ll\alpha}\frac{1}{\beta
!(\alpha-\beta)!}=\left(
\begin{array}{c}
  p+q \\
  p \\
\end{array}\right)\frac{1}{\alpha !}$ and
$\sum_{\beta'\ll\alpha'}1\leq \left(
\begin{array}{c}
  p'+q' \\
  p' \\
\end{array}\right)$. Therefore
$$ \|C\|=(\sum_{\alpha,\alpha'}\frac{\vert C_{\alpha,\alpha'}\vert^{\rho}}{\alpha
!((p+q)!(p'+q')!)^{\rho-1}})^{1/\rho} \leq $$
$$(\sum_{\alpha,\alpha'}\frac{1}{{\alpha
!((p+q)!(p'+q')!)^{\rho-1}}}(\alpha !)^{\rho}\sum_{\beta
\ll\alpha,\beta'\ll\alpha'}\frac{\vert
A_{\beta,\beta'}\vert^{\rho}}{\beta !}\frac{\vert
B_{\alpha-\beta,\alpha'-\beta'}\vert^{\rho}}{(\alpha-\beta)!}(\left(
\begin{array}{c}
  p+q \\
  p \\
\end{array}\right)\frac{1}{\alpha !}\left(
\begin{array}{c}
  p'+q' \\
  p' \\
\end{array}\right))^{\rho/\varrho})^{1/\rho} =$$
$$(\sum_{\beta \beta'}\frac{\vert
A_{\beta,\beta'}\vert^{\rho}}{\beta !(p!p'!)^{\rho
-1}})^{1/\rho}(\sum_{\gamma \gamma'}\frac{\vert B_{\gamma
\gamma'}\vert^{\rho}}{\gamma !(q!q'!)^{\rho
-1}})^{1/\rho}=\|A\|\|B\|$$ due to $\rho/\varrho =\rho-1$

To show the inequality $\lambda(p, p', q, q')\|A\|\|B\| \leq
\|A\bigodot B\|$ let us consider $$X=\{(A,B): A\in
Mat(p,p';F),\|A\|=1, B\in Mat(q,q';F),\|B\|=1\}$$, which is a
compact set
 in the corresponding finite dimensional vector space, and the  continuous map $(A,B)\mapsto A\bigodot B$.
 The image of $X$, with respect to this map, is a compact set which doesn't contain zero vector
 due to Proposition 1. Let $\lambda(p, p', q, q')> 0$ stand for the distance between zero vector and
 this image set with respect to the corresponding $\rho$- norm. So $\lambda(p, p', q, q')\leq \|A\bigodot B\|$
 for any $(A,B)\in X$ and due to Proposition 1 one has
 $$\lambda(p, p', q, q')\|A\|\|B\| \leq \|A\bigodot B\|$$ for any $A\in Mat(p,p';F), B\in Mat(q,q';F)$

{\bf Remark.} It would be nice if one could offer an expression
for $$\lambda(p, p', q,
q')=\inf_{\|A(p,p')\|=\|B(q,q')\|=1}\|A\bigodot B\|_{\rho}$$ in
terms of $n$, $n'$, $p$, $p'$, $q$, $q'$ and $\rho$.

 With respect to the ordinary product of matrices a result similar to $\|A\bigodot B\|\leq
  \|A\|\|B\|$ is not valid. But one can have the following result.

  {\bf Proposition 3.} The following inequality
   $$\|A(p,q)B(q,q')\|\leq (q!)^{2-1/\rho}(q'!)^{2/\rho -1}\|A\|\|B\|_{\varrho}$$
   is true.

 {\bf Proof.} Indeed due to the H\"{o}lder inequality one has
 $$\|A(p,q)B(q,q')\|^{\rho}=\sum_{\alpha,\alpha'}\frac{1}{\alpha
!(p!q'!)^{\rho
-1}}\vert\sum_{\beta}A_{\alpha,\beta}B_{\beta,\alpha'}\vert^{\rho}
\leq \sum_{\alpha,\alpha'}\frac{1}{\alpha !(p!q'!)^{\rho
-1}}\sum_{\beta}\vert
A_{\alpha,\beta}\vert^{\rho}(\sum_{\gamma}\vert
B_{\gamma,\alpha'}\vert^{\varrho})^{\rho/\varrho}=$$
$$\sum_{\alpha,\beta}\frac{\vert A_{\alpha,\beta}\vert^{\rho}}{\alpha
!(p!q!)^{\rho -1}}(\sum_{\gamma,\alpha'}\frac{\vert
B_{\gamma,\alpha'}\vert^{\varrho}}{\gamma !(q!q'!)^{\varrho
-1}}\gamma !)^{\rho/\varrho}(q!)^{\rho}(q'!)^{2-\rho}\leq
\|A\|^{\rho}\|B\|_{\varrho}^ {\rho}(q!)^{2\rho -1}(q'!)^{2-\rho}$$
as far as $\gamma !\leq q!$.

In particular case the following estimation is also true.

 {\bf Proposition 4.} For any nonnegative integer numbers $m, k,
q'$ and $h\in Mat_{n,n}(0,1;F)$, \\ $A\in Mat_{n,n'}(m+k,q';F)$
the following inequality
 $$\|(\frac{h^{(m)}}{m!}\bigodot E_k)A\|\leq \left(\begin{array}{c}
  m+k \\
  k \\\end{array}\right)\|h\|_{\varrho}^m\|A\|$$ is valid.

{\bf Proof.} Indeed  $$\|(\frac{h^{(m)}}{m!}\bigodot
E_k)A\|^{\rho}=\sum_{\alpha,\alpha'}\frac{1}{\alpha !(k!q'!)^{\rho
-1}}\vert((\frac{h^{(m)}}{m!}\bigodot
E_k)A)_{\alpha,\alpha'}\vert^{\rho}=\sum_{\alpha,\alpha'}\frac{1}{\alpha
!(k!q'!)^{\rho -1}}\vert\sum_{\beta}(\frac{h^{(m)}}{m!}\bigodot
E_k)_{\alpha,\beta}A_{\beta,\alpha'}\vert^{\rho}= $$

$$\sum_{\alpha,\alpha'}\frac{1}{\alpha
!(k!q'!)^{\rho
-1}}\vert\sum_{\beta}\frac{h^{\beta-\alpha}}{(\beta-\alpha)!^{1/\varrho}}
\frac{A_{\beta,\alpha'}}{(\beta-\alpha)!^{1/\rho}}\vert^{\rho}$$
as far as
$$(\frac{h^{(m)}}{m!}\bigodot
E_k)_{\alpha,\beta}=\frac{h^{\beta-\alpha}}{(\beta-\alpha)!}$$

Due to the H\"{o}lder inequality
$$(\sum_{\beta}\vert\frac{h^{\beta-\alpha}}{(\beta-\alpha)!^{1/\varrho}}
\frac{A_{\beta,\alpha'}}{(\beta-\alpha)!^{1/\rho}}\vert)^{\rho}\leq
\sum_{\beta}\frac{\vert
A_{\beta,\alpha'}\vert^{\rho}}{(\beta-\alpha)!}
(\sum_{\beta}\frac{\vert
h^{\varrho(\beta-\alpha)}\vert}{(\beta-\alpha)!})^{\rho/\varrho}
=\sum_{\beta}\frac{\vert
A_{\beta,\alpha'}\vert^{\rho}}{(\beta-\alpha)!}
(\frac{\|h\|_{\varrho}^{m\varrho}}{m!})^{\rho/\varrho}$$ Therefore
$$\|(\frac{h^{(m)}}{m!}\bigodot
E_k)A\|^{\rho}\leq \sum_{\beta,\alpha'}\frac{\vert
A_{\beta,\alpha'}\vert^{\rho}}{\beta! ((m+k)!q'!)^{\rho
-1}}\left(\begin{array}{c}
  m+k \\
  k \\\end{array}\right)^{\rho -1}\sum_{\alpha}\left(\begin{array}{c}
  \beta \\
  \alpha \\\end{array}\right) \|h\|_{\varrho}^{m\rho}=$$
  $$\sum_{\beta,\alpha'}\frac{\vert A_{\beta,\alpha'}\vert^{\rho}}{\beta!
((m+k)!q'!)^{\rho-1}}\left(\begin{array}{c}
  m+k \\
  k \\\end{array}\right)^{\rho}\|h\|_{\varrho}^{m\rho}= \|A\|^{\rho}\left(\begin{array}{c}
  m+k \\
  k \\\end{array}\right)^{\rho}\|h\|_{\varrho}^{m\rho}$$

   {\bf Corollary.} For any nonnegative integer number $m$, $q'$,  $A\in Mat_{n,n'}(m,q';F)$ and
   $h^i\in Mat_{n,n}(0,1;F)$, $i=1,2,...,m$
  the following inequality
 $$\|\frac{h^1\bigodot h^2\bigodot ... \bigodot h^m}{m!}A\|\leq \|h^1\|_{\varrho}\|h^2\|_{\varrho}
 ...\|h^m\|_{\varrho}\|A\|$$ is valid.

From now let us assume that $x_1,x_2,...,x_n$  are variables over
$F$, $q'$ is a fixed nonnegative integer and
$x=(x_1,x_2,...,x_n)\in M_{n,n}(0,1;F[x])$.

Now we are going to consider an application of the new product to
power series $\sum_{\alpha \in I_n}x^{\alpha}a_{\alpha}$, where
$a_{\alpha}\in Mat(0,q';F)$. To do it we represent the power
series $\sum_{\alpha \in I_n}x^{\alpha}a_{\alpha}$ in the form
$\sum_{m=0}^{\infty}\frac{x^{(m)}}{m!}A(m)$, where $A(m)\in
M(m,q';F)$.

It is well known that all $\rho$-norms define the same topology in
$F^n$. Therefore one can speak about convergence of the above
power series without refereing  to any particular $\rho$-norm.

{\bf Definition 3.} A power series
$\sum_{m=0}^{\infty}\frac{x^{(m)}}{m!}A(m)$ is said to be absolute
convergent at $h\in F^n$ if its each component is absolute
convergent at $h\in F^n$ e.i. for each $\alpha'\in I_{n'}$,
$|\alpha'|=q'$, the positive series $\sum_{\alpha \in
I_n}|\frac{h^{\alpha}}{\alpha !}A(m)_{\alpha, \alpha'}|$
converges.

 Due to the inequality
$$|\frac{h^{(m)}}{m!}A(m)|\leq \|h\|_{\varrho}^m\|A(m)\|$$
(Proposition 4) the power series
\begin{equation}\sum_{m=0}^\infty\frac{x^{(m)}}{m!}A(m)\end{equation} absolutely converges
whenever $\|x\|_{\varrho}< R=\frac{1}{r}$, where
$r=\overline{\lim}_{m\rightarrow\infty}\|A(m)\|_{\rho}^{\frac{1}{m}}$.

{\bf Theorem 2.}  Power series (1) is absolute convergent at $h\in
F^n$ whenever $\|h\|_{\varrho}< R$ and  for any $R_1>
Rn^{\frac{\rho -1}{\rho}}$ there exists such $\overline{h}\in F^n$
that $\|\overline{h}\|_{\varrho}= R_1$ and power series (1) is not
absolute convergent at $\overline{h}$

{\bf Proof.}

If $R=\infty$ there is nothing to prove. Assume that $1\leq \rho <
\infty$, $R < \infty$, $R_1> Rn^{\frac{\rho -1}{\rho}}$ and for
any $\overline{h}\in F^n$ for which $\|\overline{h}\|_{\varrho}=
R_1$ power series (1) is absolute convergent at $\overline{h}$.
Due to convergence of the series

 $\sum_{m=0}^{\infty}\sum_{|\alpha
|=m}|\frac{\overline{h}^{\alpha}}{\alpha!}A(m)_{\alpha,\alpha'}|$
for big enough $m$ one has
$$|\frac{\overline{h}^{\alpha}}{\alpha!}A(m)_{\alpha,\alpha'}|\geq
|\frac{\overline{h}^{\alpha}}{\alpha!}A(m)_{\alpha,\alpha'}|^{\rho}\geq
|\overline{h}^{\alpha}|^{\rho}\frac{|A(m)_{\alpha,\alpha'}|^{\rho}}{\alpha!(m!q'!)^{\rho
-1}} $$ and therefore the series $$\sum_{|\alpha'
|=q'}\sum_{m=0}^{\infty}\sum_{|\alpha
|=m}|\overline{h}^{\alpha}|^{\rho}\frac{|A(m)_{\alpha,\alpha'}|^{\rho}}{\alpha!(m!)^{\rho
-1}} $$  should converge whenever $\overline{h}\in F^n$ and
$\|\overline{h}\|_{\varrho}= R_1$. Consider this series for
$\overline{h}=(R_1n^{\frac{-1}{\varrho}},R_1n^{\frac{-1}{\varrho}},...,R_1n^{\frac{-1}{\varrho}})\in
F^n$ for which $\|\overline{h}\|_{\varrho}= R_1$.
$$\sum_{|\alpha'
|=q'}\sum_{m=0}^{\infty}\sum_{|\alpha
|=m}|\overline{h}^{\alpha}|^{\rho}\frac{|A(m)_{\alpha,\alpha'}|^{\rho}}{\alpha!(m!q'!)^{\rho
-1}}=
\sum_{m=0}^{\infty}(R_1n^{\frac{-1}{\varrho}})^{m\rho}\sum_{|\alpha
|=m,|\alpha'
|=q'}\frac{|A(m)_{\alpha,\alpha'}|^{\rho}}{\alpha!(m!q'!)^{\rho
-1}}=\sum_{m=0}^{\infty}(R_1n^{\frac{-1}{\varrho}})^{m\rho}\|A(m)\|^{\rho}
$$
 But
the radius of convergence of the ordinary number series
$\sum_{m=0}^{\infty }t^m \|A(m) \|^{\rho}$ equals $R^{\rho}$ and
$(R_1n^{\frac{-1}{\varrho}})^{\rho}>R^{\rho}$ this contradiction
indicates that the Theorem is true in this case.

Let us consider now the $\rho=\infty$ case.  Assume that $R_1>Rn$
and power series (1) is absolute convergent at any
$\overline{h}\in F^n$ for which $\|\overline{h}\|_1= R_1$. In
particular for any $\alpha'\in I_{n'}$, $|\alpha'|=q'$ and
$\overline{h}=(R_1n^{\frac{-1}{\varrho}},R_1n^{\frac{-1}{\varrho}},...,R_1n^{\frac{-1}{\varrho}})\in
F^n$, for which $\|\overline{h}\|_{1}= R_1$, the series
 $$\sum_{m=0}^{\infty}\sum_{|\alpha
|=m}|\frac{\overline{h}^{\alpha }}{\alpha!}A(m)_{\alpha,\alpha'}|=
\sum_{m=0}^{\infty} (\frac{R_1}{n})^m \sum_{|\alpha
|=m}|\frac{A(m)_{\alpha,\alpha'}}{\alpha !}|$$ is convergent. In
this case the series $$\sum_{m=0}^{\infty} (\frac{R_1}{n})^m
\sum_{|\alpha |=m,|\alpha'
|=q'}|\frac{A(m)_{\alpha,\alpha'}}{m!q'!}|$$ is convergent as
well. Due to the equality $\sum_{|\alpha
|=m}1=\left(\begin{array}{c}
  m+n-1 \\
  n-1 \\\end{array}\right)$ the following inequality
$$\|A(m)\|_{\infty}\leq \sum_{|\alpha |=m,|\alpha'
|=q'}|\frac{A(m)_{\alpha,\alpha'}}{m!q'!}|\leq
\|A(m)\|_{\infty}\left(\begin{array}{c}
  m+n-1 \\
  n-1 \\\end{array}\right)\left(\begin{array}{c}
  q'+n'-1 \\
  n'-1 \\\end{array}\right)$$ is clear. It implies that the radius
  of convergence of the power series $\sum_{m=0}^{\infty} t^m \sum_{|\alpha |=m,|\alpha'
|=q'}|\frac{A(m)_{\alpha,\alpha'}}{m!q'!}|$ is the same $R$. But
in our case $\frac{R_1}{n}> R$ and $\sum_{m=0}^{\infty}
(\frac{R_1}{n})^m \sum_{|\alpha |=m,|\alpha'
|=q'}|\frac{A(m)_{\alpha,\alpha'}}{m!q'!}|$ converges. This
contradiction indicates that the Theorem is true in this case as
well.

 Due to this result if for power series (1) the radius of convergence $R$ is zero then in any neighborhood
 of zero one can find a point where (1) is not absolute convergent. This theorem indicates also a privileged
  position of $1$-norm among all
$\rho$-norms as  far as in this case $Rn^{\frac{\rho -1}{\rho}}=R$
and for $\rho >1$ power series (1) has a solid layer  $$\{h\in
F^n: R\leq \|h\|_\varrho \leq Rn^{\frac{\rho -1}{\rho}}\}$$ of
indeterminacy. Investigation the behavior of power series (1)in
this layer of indeterminacy could be interesting. It is hoped that
the emergence of a layer of indeterminacy is not a defect of our
approach.

{\bf Question 1.} Let $\sum_{\alpha}a_{\alpha}x^{\alpha}$ be any
series for which $\sum_{m=0}^{\infty}H_m(x)$, where
$H_m(x)=\sum_{|\alpha |=m }a_{\alpha}x^{\alpha}$, is absolute
convergent in some neighborhood of zero. Does it imply that the
original series $\sum_{\alpha}a_{\alpha}x^{\alpha}$ is also
absolute convergent in some neighborhood of zero as well?

{\bf Question 2.} Consider power series (1), for each nonnegative
$m$ multilinear map $$(F^n)^m \longrightarrow F^{n'}:
 (h^1, h^2, ... , h^m)\longmapsto (h^1\bigodot h^2\bigodot ... \bigodot h^m)A(m)$$ and
 its norm $\|A(m)\|$ ([2]) defined by
$$\|A(m)\|= \sup_{\|h^1\|_\varrho =\|h^2\|_\varrho =...=\|h^m\|_\varrho =1}
\|(h^1\bigodot h^2\bigodot ... \bigodot h^m)A(m)\|_\rho$$
According to the above Corollary for each $m$ the inequality
$\|A(m)\|\leq  \|A(m)\|_\rho$ is true. Is it true that
$$\overline{\lim}_{m\rightarrow\infty}\|A(m)\|^{\frac{1}{m}}=\overline{\lim}_{m\rightarrow\infty}\|A(m)\|_{\rho}^{\frac{1}{m}}$$

\begin{center}{References}\end{center}

[1] U.Bekbaev. A matrix representation of composition of
polynomial maps. arXiv:math 0901.3179v3

[2] H.Cartan. Calcul diff\'{e}rntiel. Forms diff\'{e}rentielles.
Hermann. Paris, 1967.

 \end{document}